\documentclass{amsart}

\usepackage{hyperref}
\usepackage{amsmath, amsthm, amssymb}
\usepackage{url}
\usepackage{graphicx}

\newtheorem{theorem}{Theorem}
\newtheorem{question}{Question}

\theoremstyle{definition}
\newtheorem*{remark}{Remark}

\begin{document}

\title{On Functions Whose Mean Value Abscissas Are Midpoints, with Connections to Harmonic Functions}%

\author{Paul Carter and David Lowry-Duda}

\maketitle

\begin{abstract}
  We investigate functions with the property that for every interval, the slope at the midpoint of the interval is the same as the average slope.
  More generally, we find functions whose average slopes over intervals are given by the slope at a weighted average of the endpoints of those intervals.
  This is equivalent to finding functions satisfying a weighted mean value property.
  In the course of our exploration, we find connections to harmonic functions that prompt us to explore multivariable analogues and the existence of ``weighted harmonic functions.''
\end{abstract}

%%% Begin document %%%

\section{Introduction and Statement of the Problem}

Recall the mean value theorem of calculus, which states that if $f$ is a differentiable function on $[a,b]$, then there is a $c \in (a,b)$, which we call the \emph{mean value abscissa}, such that
\begin{equation*}
  \frac{f(b) - f(a)}{b-a} = f'(c).
\end{equation*}
In typical applications, such as in most proofs of the fundamental theorems of calculus, we are given $f$ and know there exists such a $c$.
But what if we place conditions on $c$ and ask when there exist functions $f$ satisfying the mean value theorem with that $c$?
For instance, what if $c$ is exactly halfway between $a$ and $b$?

\begin{question}\label{q1}
  For which differentiable functions $f$ is
\begin{equation}\label{eq:MAIN}
  \frac{f(b) - f(a)}{b - a} = f'\left( \frac{b+a}{2} \right)
\end{equation}
for all $b > a$?
\end{question}
\noindent In other words, for which functions $f$ does the midpoint of each interval serve as the mean value abscissa for that interval?
The midpoint of $[a,b]$ is the average of $a$ and $b$.
What about more general averages, like weighted averages?

\begin{question}\label{q2}
  For a fixed $\lambda \in (0,1)$, which differentiable functions $f$ satisfy
\begin{equation}\label{eq:MAIN2}
  \frac{f(b) - f(a)}{b - a} = f' \left(\lambda a + (1 - \lambda)b \right)
\end{equation}
for all $b > a$?
\end{question}
\noindent In other words, for which functions $f$ do weighted averages of the endpoints of intervals serve as mean value abscissas for those intervals?

We explore the answers to these two questions.
Along the way, we'll find strong connections between functions satisfying these constrained mean value theorem requirements and harmonic functions which will prompt us to investigate high dimensional analogues.

\section{A First Attempt at the First Question}

Before trying to find all such functions, we should note that if $f(x)$ is an affine function $f(x) = ax + b$ for any $a,b$, then $f$ satisfies the conditions for both questions.
So there are infinite families of solutions.
But without further details, it might be hard to find all functions satisfying~\eqref{eq:MAIN}.

\begin{remark}
  In fact, it was known to Archimedes that parabolas also satisfy~\eqref{eq:MAIN}.
  This appears in Proposition~1 in his work on the quadrature of the parabola~\cite{heath1897works}.
  In modern times, some calculus textbook authors include verifying that parabolas satisfy~\eqref{eq:MAIN} as an exercise.
  For instance, it appears as an ``Additional and Advanced'' exercise in~\cite{thomas2010thomas}. %Chapter 3 Review Problem 14
\end{remark}

When first encountering a problem, it can be useful to consider simplifying assumptions.
We might be able to find more functions by expanding the left- and right-hand sides of~\eqref{eq:MAIN}.
Suppose that $f$ is a function that satisfies~\eqref{eq:MAIN}, and suppose $f$ is at least $3$ times continuously differentiable.
Then we have the Taylor polynomial $f(x) = f(y) + f'(y)(x-y) + \tfrac{1}{2}f''(y)(x-y)^2 + \tfrac{1}{6}f'''(y)(x-y)^3 + h(x)(x-y)^3$, where $h(x) \to 0$ as $x \to y$.
Subtracting $f(y)$ and dividing by $x-y$ yields an expression for the left-hand side of~\eqref{eq:MAIN}:
\begin{equation*}
  \frac{f(x)- f(y)}{x-y} = f'(y) + \tfrac{1}{2}f''(y)(x-y) + \tfrac{1}{6}f'''(y)(x-y)^2 + h(x)(x-y)^2.
\end{equation*}
Taking the Taylor polynomial for $f'(x)$ centered at $y$ and substituting $(x+y)/2$, we see that the right-hand side of~\eqref{eq:MAIN} is
\begin{equation*}
  f'(y) + \tfrac{1}{2} f''(y) (x-y) + \tfrac{1}{8} f'''(y)(x-y)^2 + \tfrac{1}{4}\widetilde{h}(\tfrac{x+y}{2})(x-y)^2,
\end{equation*}
where $\widetilde{h}(x) \to 0$ as $x \to y$.
Setting these two expressions equal, we cancel the $f'(y)$ and $f''(y)$ terms and divide each side by $(x-y)^2$.
Letting $x \to y$, we see that $\tfrac{1}{6}f'''(y) = \tfrac{1}{8} f'''(y)$, or rather that $f'''(y) = 0$.

Since this is true for all $y$, the third derivative is identically $0$.
So $f$ is at most a quadratic polynomial. We can check that if $f(x) = a + bx + cx^2$, then
\begin{equation*}
  \frac{f(x) - f(y)}{x-y} = \frac{b(x-y) + c(x^2 - y^2)}{x-y} = b + c(x+y) = f'\left(\frac{x+y}{2}\right),
\end{equation*}
so every polynomial of degree at most $2$ is an answer to Question 1.

Exploration through Taylor polynomials led us only to parabolas.
But we made a big assumption.
What about $f$ that satisfy the conditions of the first question but which are not necessarily three times differentiable?

\subsection{Connection to harmonic functions.}

By fixing various $a$ and allowing $b$ to vary in~\eqref{eq:MAIN}, one can show that $f'(x)$ is continuous.
We now use the fundamental theorem of calculus to rewrite~\eqref{eq:MAIN} as
\begin{equation}
  \frac{1}{b-a} \int_a^b f'(t) \, dt = f' \left(\frac{b+a}{2}\right),
\end{equation}
or equivalently as
\begin{equation}\label{eq:fderiv_is_harmonic}
  \frac{1}{2h} \int_{x-h}^{x+h} f'(t)\, dt = f'(x).
\end{equation}
This can be interpreted to mean that the value $f'(x)$ is given by the average of $f'$ in any symmetric interval around $x$, a very restrictive property.
This is precisely the mean value property satisfied by harmonic functions,
\begin{equation}\label{eq:harmonic_mvt}
  g(x) = \frac{1}{|B_h(x)|} \int_{B_h(x)} g(t)\, dV,
\end{equation}
where $g$ is a harmonic function on $\mathbb{R}^n$, $B_h(x)$ denotes the ball of radius $h$ centered at $x$, $dV$ denotes the standard Euclidean volume measure, and $|B_h(x)|$ denotes the volume of $B_h(x)$.
Harmonic functions appear in many areas of mathematics, but are particularly fundamental to the theory of differential equations.
Recall that a twice differentiable function $g$ is called \emph{harmonic} if the sum of the second partial derivatives of $g$ is identically $0$.

At least, that's the standard definition.
Harmonic functions can be thought of in more intuitive ways.
One particularly illuminating overview on harmonic functions is given by Needham~\cite{needham1994geometry}, who explores geometric interpretations of harmonic functions by exploring relations similar to~\eqref{eq:harmonic_mvt}.
In the memorably titled~\cite{kac1966can}, Kac gives a different approach to harmonic functions while explaining why they are called ``harmonic.''

The study of harmonic functions usually begins with $\mathbb{R}^2$, as it is trivial to identify the harmonic functions on $\mathbb{R}$.
The only functions on $\mathbb{R}$ with identically zero second derivative are affine functions $ax + b$.
It is a remarkable theorem that any function satisfying~\eqref{eq:harmonic_mvt} is itself harmonic (for the standard proof, see Lemma~4.6 and subsequent discussion in~\cite{stein2009real}).
In fact, the condition~\eqref{eq:MAIN} is the degenerate one-dimensional form of another integral mean value property shared by harmonic functions, namely
\begin{equation}\label{eq:meanvalueshell}
  g(x) = \frac{1}{|\partial B_h(x)|} \int_{\partial B_h(x)} g(t) \, dS,
\end{equation}
where $\partial B_h(x)$ denotes the boundary of $B_h(x)$, $|\partial B_h(x)|$ is the surface area of $B_h(x)$, and $dS$ denotes the standard Euclidean surface area measure.

So for any function $f$ satisfying~\eqref{eq:MAIN}, its derivative $f'$ satisfies~\eqref{eq:fderiv_is_harmonic}, and so is harmonic.
Thus $f'$ is at most a linear polynomial, so that $f$ is at most a quadratic polynomial.
We have filled in the gap left from Taylor expansions by appealing to harmonic functions.

\section{Investigating Weighted Mean Values}

Let us now try to understand the more general Question~\ref{q2}.
Pursuing the connection to harmonic functions, we can rewrite~\eqref{eq:MAIN2} as the integral weighted mean value property
\begin{equation}\label{eq:weighted_mean_integral_badcoords}
  f'\left(\lambda a + (1-\lambda) b\right) = \frac{1}{b-a} \int_a^b f'(t)\, dt.
\end{equation}
This is similar to the regular harmonic mean value property~\eqref{eq:fderiv_is_harmonic}, and so we might ask whether functions satisfying~\eqref{eq:weighted_mean_integral_badcoords} are harmonic.
This turns out to be the case.

To prove this, we first show that any function satisfying~\eqref{eq:weighted_mean_integral_badcoords} is infinitely differentiable.
We take $x=\frac{a+b}{2}$ and $h=\frac{b-a}{2}$ and rewrite~\eqref{eq:weighted_mean_integral_badcoords} as
\begin{equation}\label{eq:weighted_mean_integral_goodcoords}
    f'\left(x+(1-2\lambda) h\right)=\frac{1}{2h} \int_{x-h}^{x+h} f'(t)\, dt,
\end{equation}
which now holds for any $x$ and any $h > 0$.
We can center the left hand side at $x$ through a change of variables, getting
\begin{equation*}
  f'(x) = \frac{1}{2h} \int_{x - (2 - 2\lambda)h}^{x + 2\lambda h} f'(t)\, dt.
\end{equation*}
By the fundamental theorem of calculus, we can directly differentiate $f'(x)$ to see that $f''(x) = \frac{1}{2h} \left( f'(x + 2\lambda h) - f'(x - (2 - 2\lambda)h )\right)$.
By the same argument, each of the two terms appearing in the derivative are differentiable with derivatives expressible as linear combinations of $f'$.
Inductively, we can show that $f'$ is infinitely differentiable.

\begin{remark}
  We can see the infinite differentiability in a different way.
  We can compose the integral weighted mean value property with itself to write $f'$ as the $2$-fold integral
  \begin{equation*}
    f'(x) = \frac{1}{2h} \int_{x - (2 - 2\lambda)h}^{x + 2\lambda h} \left( \frac{1}{2h} \int_{t - (2 - 2\lambda h)}^{t + 2\lambda h} f'(u) \, du \right) \, dt,
  \end{equation*}
  which is clearly twice-differentiable, again by the fundamental theorem of calculus.
  Composing $n$ times, we can represent $f'$ as an $n$-fold integral which is $n$ times differentiable.
  So $f'$, and therefore $f$, is infinitely differentiable.
\end{remark}

As $f$ is at least twice continuously differentiable, we can use Taylor's Theorem as before to write $f(x) = f(y) + f'(y)(x-y) + \tfrac{1}{2}f''(y)(x-y)^2 + h(x)(x-y)^2$, where $h(x) \to 0$ as $x \to y$.
Similar to when we first tackled Question~1, we expand the left-hand side and right-hand side of~\eqref{eq:MAIN2}.
The left-hand side remains
\begin{equation*}
  f'(y) + \tfrac{1}{2}f''(y)(x-y) + h(x)(x-y).
\end{equation*}
Using the linear Taylor polynomial for $f'(x)$ and simplifying, the right-hand side becomes
\begin{equation*}
  f'(y) + \lambda f''(y)(x-y) + \lambda \widetilde{h}(\lambda x + (1-\lambda)y)(x-y),
\end{equation*}
where $\widetilde{h}(x) \to 0$ as $x \to y$.
Setting these equal, we may cancel $f'(y)$ and divide by $(x-y)$.
Letting $x \to y$, we see that $\tfrac{1}{2}f''(y) = \lambda f''(y)$.

There are two possibilities.
If $\lambda = \tfrac{1}{2}$, then our ``weighted average'' is just the normal average from Question 1.
Otherwise, we must have $f''(y) = 0$.
This is true for all $y$, and so the second derivative is identically $0$.
We conclude that the only twice differentiable functions satisfying the conditions from Question~2 are linear polynomials $f(x) = ax + b$ unless $\lambda = \tfrac{1}{2}$, when $f$ can be a quadratic polynomial.

We have shown that the normally-weighted mean value property on $\mathbb{R}$ is distinguished as the only weighted mean value property that leads to a nontrivial family.
It was not obvious (to the authors, at least) that there would be no analogous ``weighted harmonic functions'' corresponding to a weighted mean value property, and this adds to the list of special properties held by harmonic functions.

We summarize the answers to the two original questions in the following theorem.

\begin{theorem}\label{thm:firsttwoquestions}
  Fix a $\lambda \in (0,1)$, and suppose $f:\mathbb{R} \longrightarrow \mathbb{R}$ satisfies the weighted mean value condition
  \begin{equation}
    \frac{f(b) - f(a)}{b-a} = f'(\lambda a + (1-\lambda) b)
  \end{equation}
  for all $a < b$. Then $f$ is a quadratic polynomial. Further, if $\lambda \neq \tfrac{1}{2}$, then $f$ is merely a linear polynomial.
\end{theorem}

\begin{remark}
  We have used that~\eqref{eq:MAIN2} is supposed to hold for all $a < b$, allowing both $a$ and $b$ to vary.
  Another interesting question comes from fixing $a=0$ and letting $b$ vary.
  Then we might look for $f$ with the mean value abscissa on the intervals $[0, b]$ being given by weighted averages of $0$ and $b$.
  In this case, there \emph{are} additional families of functions, but only for $\lambda = 1/(k+1)^{1/k}$.
  The reader might try to explore this approach using Taylor expansions, or perhaps another method entirely.
\end{remark}

\section{Weighted Harmonic Functions in Higher Dimensions}

Our initial questions led us to study functions satisfying the integral weighted mean value property for $\mathbb{R}$, and we found that the non-weighted average has a distinguished connection to harmonic functions.
Since there are many more harmonic functions on $\mathbb{R}^n$ for $n \geq 2$ than there are for $\mathbb{R}$, there might be enough wiggle room for a ``weighted harmonic function'' to exist corresponding to a weighted mean value property on $\mathbb{R}^n$.
In analogy with~\eqref{eq:weighted_mean_integral_goodcoords}, we can ask about the following integral weighted mean value property
\begin{equation}\label{eq:weighted_mean_Rn}
  g(\mathbf{x} + (1 - 2\lambda)h\mathbf{v}) = \frac{1}{|B_h(\mathbf{x})|} \int_{B_h(\mathbf{x})} g(\mathbf{t}) \, dV.
\end{equation}
Here, $\lambda \in (0,1)$ is the weight, and $\mathbf{v}$ is a unit vector indicating the direction of the weighting.
The mean value abscissa $(\mathbf{x} + (1 - 2\lambda)h\mathbf{v})$ is a point in the ball $B_h(\mathbf{x})$ that differs from the center of the ball by a distance proportional to the radius.
Then~\eqref{eq:weighted_mean_Rn} can be interpreted to mean that the average value of $g$ on the ball $B_h(\mathbf{x})$ is given by the value of $g$ at the point $(\mathbf{x} + (1 - 2\lambda)h\mathbf{v})$ inside that ball.
So we ask the following question.

\begin{question}\label{q3}
  For a fixed $\lambda \in (0,1)$ and a fixed unit vector $\mathbf{v} \in \mathbb{R}^n$, which functions $g \in C^1(\mathbb{R}^n)$ satisfy
  \begin{equation*}
    g(\mathbf{x} + (1 - 2\lambda)h\mathbf{v}) = \frac{1}{|B_h(\mathbf{x})|} \int_{B_h(\mathbf{x})} g(\mathbf{t}) \, dV
  \end{equation*}
  for all $\mathbf{x}$ and all $h > 0$?
\end{question}

Notice that when $n = 1$ and $v = 1$, this is exactly Question~\ref{q2}.
As this generalizes Question~\ref{q2}, our method also generalizes our approach to answering Question~\ref{q2}.

When $\lambda = \frac{1}{2}$, the mean value abscissa is exactly the center of the ball, and so the weighted mean value property becomes the ordinary mean value property~\eqref{eq:harmonic_mvt}.
Then for $\lambda = \frac{1}{2}$, it is exactly harmonic functions $g$ which answer Question~\ref{q3}.

Let us now consider $\lambda \neq \frac{1}{2}$.
We first show that any function satisfying~\eqref{eq:weighted_mean_Rn} is infinitely differentiable.
Rewrite~\eqref{eq:weighted_mean_Rn} as
\begin{equation}\label{eq:weighted_mean_Rn_shifted}
  g(\mathbf{x}) = \frac{1}{|B_h|} \int_{B_h(\mathbf{x}-(1 - 2\lambda)h \mathbf{v})} g(\mathbf{t}) \, dV,
\end{equation}
where $\lvert B_h \rvert$ denotes the volume of the ball of integration.
As $g$ is continuously differentiable, we can now compute any partial derivative of $g$ by differentiating~\eqref{eq:weighted_mean_Rn_shifted} using the Leibniz rule for higher dimensions (sometimes called the Reynolds transport theorem, see~\cite{flanders1973} for additional exposition) to obtain
\begin{align*}
  \partial_ig(\mathbf{x}) &= \frac{1}{|B_h|} \int_{\partial B_h(\mathbf{x}-(1 - 2\lambda)h \mathbf{v})} g(\mathbf{t}) \mathbf{e}_i\cdot \mathbf{n} \, dS \\
  &=\frac{1}{|B_h|} \int_{\partial B_h(0)} g(\mathbf{t}+\mathbf{x}-(1 - 2\lambda)h \mathbf{v}) \mathbf{e}_i\cdot \mathbf{n} \, dS,
\end{align*}
where $\mathbf{e}_i$ denotes the $i$-th standard basis vector of $\mathbb{R}^n$, $\partial_i$ denotes the partial derivative with respect to the $i$-th coordinate, and $\mathbf{n}$ is the outward-pointing unit normal vector field of the boundary $\partial B_h$.
To compute the second partial derivative $\partial_{ij}$, we can now differentiate under the integral sign
\begin{align*}
  \partial_{ij}g(\mathbf{x}) &= \frac{1}{|B_h|} \int_{\partial B_h(0)} \partial_jg(\mathbf{t}+\mathbf{x}-(1 - 2\lambda)h \mathbf{v}) \mathbf{e}_i\cdot \mathbf{n}\, dS
\end{align*}
since the first partial derivatives of $g$ exist and are continuous.
Continuing inductively, we conclude that $g$ is smooth.

  \begin{remark}
    In fact, it should be possible to relax the assumption that $g\in C^1$ and deduce this from the argument above under the assumption that $g$ is at least continuous. However, for simplicity in applying the Reynolds transport theorem, we assume $g$ has continuous partial derivatives.
  \end{remark}

As $g$ is at least twice continuously differentiable, we can again use Taylor's Theorem, which states
\begin{equation}
  g(\mathbf{x}) = g(\mathbf{y}) + (\mathbf{x}-\mathbf{y})\cdot \nabla g(\mathbf{y}) +\mathcal{O}(|\mathbf{x}-\mathbf{y}|^2),
\end{equation}
where $\nabla g$ denotes the gradient of $g$ and $\mathcal{O}$ is big-oh notation, indicating roughly that the remainder vanishes at least as quickly as $\lvert \mathbf{x}-\mathbf{y} \rvert^2$ as $\mathbf{x} \to \mathbf{y}$.
We now compute
\begin{align*}
  g(\mathbf{x}+(1-2\lambda)h \mathbf{v})&=  \frac{1}{|B_h|} \int_{B_h(\mathbf{x})} g(\mathbf{t}) \, dV \\
  &= \frac{1}{|B_h|} \int_{B_h(0)} g(\mathbf{x}+\mathbf{t}) \, dV_\mathbf{t} \\
  &= \frac{1}{|B_h|} \int_{B_h(0)} \bigg( g(\mathbf{x})+\mathbf{t}\cdot \nabla f(\mathbf{x}) +\mathcal{O}(|\mathbf{t}|^2) \bigg) \, dV_\mathbf{t} \\
  &= g(\mathbf{x})+\nabla g(\mathbf{x})\cdot\left(\frac{1}{|B_h|} \int_{B_h(0)} \mathbf{t}\, dV_\mathbf{t}\right)+\mathcal{O}(h^2).
\end{align*}
We use the notation $dV_\mathbf{t}$ to remind ourselves that the variable of integration is $\mathbf{t}$, not $\mathbf{x}$.
The integral of the vector field $\mathbf{t}$ over the ball $B_h(0)$ is zero, so the integral term above vanishes.
We have thus shown that
\begin{equation*}
  g(\mathbf{x} + (1 - 2\lambda)h \mathbf{v}) = g(\mathbf{x}) + \mathcal{O}(h^2).
\end{equation*}
Using this expression, we can show that $g$ is constant in the $\mathbf{v}$-direction, that is, $\frac{\partial g}{\partial \mathbf{v}}=0$.
We compute
\begin{equation*}
  \frac{\partial g}{\partial \mathbf{v}}=\lim_{h\to 0}\frac{g(\mathbf{x}+(1-2\lambda)h\mathbf{v})-g(\mathbf{x})}{(1-2\lambda)h} = \lim_{h\to 0}\frac{g(\mathbf{x})+\mathcal{O}(h^2)-g(\mathbf{x})}{(1-2\lambda)h} =0.
\end{equation*}
As $g$ is constant in the $\mathbf{v}$-direction, and the mean value abscissa $(\mathbf{x} + (1-2\lambda)h\mathbf{v})$ deviates from the center of the ball $B_h(\mathbf{x})$ in exactly the $\mathbf{v}$-direction, we can show that $g$ is harmonic.
At each $\mathbf{x}$, the function $g$ satisfies
\begin{equation*}
  g(\mathbf{x})=g(\mathbf{x}+(1-2\lambda)h\mathbf{v}) =  \frac{1}{|B_h(\mathbf{x})|} \int_{B_h(\mathbf{x})} g(\mathbf{t}) \, dV,
\end{equation*}
which is precisely the mean value property for harmonic functions.

By choosing an orthonormal basis for $\mathbb{R}^n$ with $\mathbf{v}$ as a basis vector, we can think of $g$ as coming from a function $\widetilde{g}$ on $\mathbb{R}^{n-1}$, which is extended to $\mathbb{R}^n$ trivially by having no dependence on the $\mathbf{v}$-coordinate.
Noting that the sums of the second partial derivatives of $g$ and $\widetilde{g}$ are the same, we have that $\widetilde{g}$ is also harmonic.
So a function $g$ satisfying the integral weighted mean value theorem~\eqref{eq:weighted_mean_Rn} for $\lambda \neq \frac{1}{2}$ is really a harmonic function on $\mathbb{R}^{n-1}$, extended to $\mathbb{R}^n$ by being constant in the $\mathbf{v}$-direction.

For example, in the case $n=2$, without loss of generality we may take $\mathbf{v}=(0,1)$, i.e., the unit vector in the $y$-direction, so that $g$ is constant in $y$.
From this we see that the weighted harmonic functions are simply those which are harmonic functions of the single variable $x$, that is, linear functions $g(x,y)=ax+b$.

In conclusion, we have shown the following
\begin{theorem}\label{thm:thirdquestion}
  Fix $\lambda \in (0,1)$ and a unit vector $\mathbf{v} \in \mathbb{R}^n$. Suppose $g: \mathbb{R}^n \longrightarrow \mathbb{R}$ satisfies
  \begin{equation*}
    g(\mathbf{x}+(1-2\lambda)h\mathbf{v}) = \frac{1}{|B_h(\mathbf{x})|} \int_{B_h(\mathbf{x})} g(\mathbf{t})\, dV
  \end{equation*}
  for all $\mathbf{x} \in \mathbb{R}^n$ and $h > 0$. Then $g$ is a harmonic function. Further, if $\lambda \neq \tfrac{1}{2}$, then $g$ is constant in the $\mathbf{v}$-direction, i.e., $\frac{\partial}{\partial \mathbf{v}} g = 0$.
\end{theorem}

This is consistent with our answers to Questions~\ref{q1} and~\ref{q2}, and includes Theorem~\ref{thm:firsttwoquestions} as the $\mathbb{R}^1$ case.
This is also consistent with our earlier observation that the normally-weighted mean value property is distinguished among weighted mean value properties by including a strictly larger family of functions.

\begin{remark}
 As a final note, it is also very natural to think of those $g$ satisfying the alternate weighted mean value property
\begin{align}\label{eq:meanvalueshellweighted}
  g(\mathbf{x} + (1 - 2\lambda)h\mathbf{v}) = \frac{1}{|\partial B_h(\mathbf{x})|} \int_{\partial B_h(\mathbf{x})} g(\mathbf{t}) \, dS,
  \end{align}
  analogous to the harmonic mean value property~\eqref{eq:meanvalueshell}. We claim that, similar to the case of harmonic functions, this condition is equivalent to~\eqref{eq:weighted_mean_Rn} and hence leads to a result equivalent to Theorem~\ref{thm:thirdquestion}. Again assuming for simplicity that $g\in C^1$, using a very similar argument as above, it is not hard to show that $g$ is in fact infinitely differentiable and satisfies $\frac{\partial g}{\partial \mathbf{v}}=0$. We do not repeat the details here, but just note that to show differentiability, we rewrite the property~\eqref{eq:meanvalueshellweighted} as
  \begin{equation}
    g(\mathbf{x}) = \frac{1}{|\partial B_h|} \int_{\partial B_h(0)} g(\mathbf{t}+\mathbf{x}-(1 - 2\lambda)h\mathbf{v}) \, dV.
  \end{equation}
  We can now differentiate this expression to obtain smoothness of $g$. From this, we use a similar argument involving Taylor's theorem as above to find that $\frac{\partial g}{\partial \mathbf{v}}=0$. From this we deduce the harmonic mean value property~\eqref{eq:meanvalueshell}, and hence the conclusion of Theorem~\ref{thm:thirdquestion} remains valid.
  \end{remark}

\section*{Acknowledgements}

D. L.-D.\ was supported by the  National Science Foundation Graduate Research Fellowship Program under Grant No. DGE 0228243.
P. C.\ gratefully acknowledges support by the National Science Foundation through grant DMS-1148284.

%\begin{biog}
%  \item[Paul Carter] obtained his PhD from Brown University in 2016. He recently moved to Tucson, where he joined the University of Arizona as a Postdoctoral Research Associate. He enjoys the good weather while rock climbing or working on differential equations.
%  \begin{affil}
%    Department of Mathematics, Brown University, Providence, RI 02906 \\ pacarter@math.brown.edu
%  \end{affil}
%
%  \item[David Lowry-Duda] earned his MS from Brown University and plans on completing his PhD there soon. He recently helped a close friend and collaborator move to Tucson. Although he normally thinks about analytic number theory and algebraic geometry, the mean value theorem seems to pop up everywhere.
%  \begin{affil}
%    Department of Mathematics, Brown University, Providence, RI 02906 \\ djlowry@math.brown.edu
%  \end{affil}
%\end{biog}
%

\begin{thebibliography}{9}

  \bibitem{evans2010partial}
    L. C. Evans,
    \emph{Partial Differential Equations}.
    Second edition,
    Graduate Studies in Mathematics,
    Vol. 19,
    American Mathematical Society,
    Providence, RI,
    2010,
    \url{http://dx.doi.org/10.1090/gsm/019}.

  \bibitem{flanders1973}
    H. Flanders,
    Differentiation under the integral sign,
    \emph{Amer. Math. Monthly} \textbf{80} no. 6 (1973) 615--627,
    \url{http://dx.doi.org/10.2307/2319163}.

  \bibitem{heath1897works}
    T. L. Heath,
    \emph{The Works of Archimedes},
    Cambridge Univ. Press,
    Cambridge,
    1897.

  \bibitem{kac1966can}
    M. Kac,
    Can one hear the shape of a drum? %, %including the comma looks silly
    \emph{Amer. Math. Monthly} \textbf{73} no. 4 (1966) 1--23,
    \url{http://dx.doi.org/10.2307/2313748}.

  \bibitem{needham1994geometry}
    T. Needham,
    The geometry of harmonic functions,
    \emph{Math. Mag.} \textbf{67} no. 2 (1994) 92--108,
    \url{http://dx.doi.org/10.2307/2690683}.

  \bibitem{stein2009real}
    E. M. Stein, R. Shakarchi,
    \emph{Real Analysis: Measure Theory, Integration, and Hilbert Spaces}.
    Princeton Lectures in Analysis,
    Princeton Univ. Press,
    Princeton, NJ,
    2009.

  \bibitem{stein2011fourier}
    %E. M. Stein, R. Shakarchi,
    ---,
    \emph{Fourier Analysis: An Introduction}.
    Princeton Lectures in Analysis,
    Princeton Univ. Press,
    Princeton, NJ,
    2011.

  \bibitem{thomas2010thomas}
    G. B. Thomas, M. D. Weird, J. Hass, F. R. Giordano,
    \emph{Thomas' Calculus Early Transcendentals}.
    Twelfth edition.
    Pearson Addison-Wesley,
    Boston, MA,
    2010.

\end{thebibliography}
\end{document}